\documentclass [winedt,yap]{iitparcgr}
\usepackage{cite}
\usepackage{amsmath,amssymb,amsfonts,bm}
\usepackage{ulem}
\usepackage[pdftex]{graphicx}       

                \usepackage{accents}
                \DeclareMathSymbol{\widehatsym}{\mathord}{largesymbols}{"62}
                \newcommand\lowerwidehatsym{%
                  \text{\smash{\raisebox{-1.3ex}{%
                    $\widehatsym$}}}}
                \newcommand\fhat[1]{%
                  \mathchoice
                    {\accentset{\displaystyle\lowerwidehatsym}{#1}}
                    {\accentset{\textstyle\lowerwidehatsym}{#1}}
                    {\accentset{\scriptstyle\lowerwidehatsym}{#1}}
                    {\accentset{\scriptscriptstyle\lowerwidehatsym}{#1}}
                }

\begin{document}\normalem
%
\frontmatter          

\IssuePrice{25.00}%
\TransYearOfIssue{2018}%
\TransCopyrightYear{2018}%
\OrigYearOfIssue{2018}%
\OrigCopyrightYear{2018}%

\TransVolumeNo{79}%
\TransIssueNo{11}
\OrigIssueNo{11}%
\OrigPages{127--149} 
\OrigCopyrightedAuthors{%
{P.Yu.~Chebotarev},
{Ya.Yu.~Tsodikova},
{A.K.~Loginov},
{Z.M.~Lezina},
{V.A.~Afonkin},
{V.A.~Malyshev}%
}
\OrigJournalName{Avtomatika i Telemekhanika}

\mainmatter

\setcounter{page}{1} 

\Rubrika{CONTROL IN SOCIAL ECONOMIC SYSTEMS}


\def\x#1{} 
\def\xz{\hspace{-.17em}}
\def\MP{M\xz P}
\def\ell{n}
\def\a{\alpha}
\def\Up#1{\vspace{-#1em}}
\def\cdc{,\ldots,}
\def\xz{\hspace{-.07em}}
\def\xy{\hspace{.07em}}
\def\beq#1{\begin{equation}\label{#1}}
\def\eeq{\end{equation}}
\def\M{\fhat{\rm M}}                             
\def\Up#1{\vspace{-#1em}}

\title{Comparative Efficiency of Altruism and Egoism as Voting Strategies in Stochastic Environment}
\author{%
{P.~Yu.~Chebotarev}$^{\xy*,**,***a}$, 
{Ya.~Yu.~Tsodikova}$^{\xy*,b}$,
{A.~K.~Loginov}$^{\xy*,c}$,
{Z.~M.~Lezina}$^{\xy*,d}$,
{V.~A.~Afonkin}$^{\xy*,***,e}$,
{V.~A.~Malyshev}$^{\xy*,**,***,f}$\\
}
\institute{$^*$Trapeznikov Institute of Control Sciences, Russian Academy of Sciences, Moscow, Russia\\
           $^{**}$Kotelnikov Institute of Radioengineering and Electronics, Russian Academy of Sciences, Moscow, Russia\\
           $^{***}$Moscow Institute of Physics and Technology, Moscow, Russia\\
           e-mail: $^a$pavel4e@gmail.com, $^b$codikova@mail.ru, $^c$a\_k\_log@mail.ru, $^d$lezinazo@gmail.com, $^e$afonkinvadim@yandex.ru, $^f$vit312@gmail.com}

\received{Received December 13, 2017}
\authorrunning{CHEBOTAREV, TSODIKOVA, LOGINOV, LEZINA, AFONKIN, MALYSHEV}
\titlerunning{EFFICIENCY OF ALTRUISM AND EGOISM AS VOTING STRATEGIES}

\maketitle

\begin{abstract}
In this paper, we study the efficiency of egoistic and altruistic strategies within the model of social dynamics determined by voting in a stochastic environment (the ViSE model) using two criteria: maximizing the average capital increment and minimizing the number of bankrupt participants. 
The proposals are generated stochastically; three families of the corresponding distributions are considered: normal distributions, symmetrized Pareto distributions, and Student's $t$-distributions.
It is found that the ``pit of losses'' paradox described earlier does not occur in the case of heavy-tailed distributions.
The egoistic strategy better protects agents from extinction in aggressive environments than the altruistic ones, however, the efficiency of altruism is higher in more favorable environments. A comparison of altruistic strategies with each other shows that in aggressive environments, everyone should be supported to minimize extinction, while under more favorable conditions, it is more efficient to support the weakest participants.
Studying the dynamics of participants' capitals we identify situations where the two considered criteria contradict each other.
At the next stage of the study, combined voting strategies and societies involving participants with selfish and altruistic strategies will be explored.


\medskip\noindent$\!\!\!\!\!\!\!\!$\DOI{XXXX} 
\end{abstract}

\setcounter{footnote}{0}

\section{Introduction}
\label{s_intro}

The ViSE (Voting in Stochastic Environment) model was proposed in 2003 and examined in \cite{1,2,3,4,5,6,7,8} and several other papers. In this study, the ViSE model is investigated in the case where participants may have an altruistic strategy aimed at supporting the poorest strata of society. The ``support window'' determines the part of society that they support; this parameter can be optimized with respect to strategy efficiency. Important external parameters are the favorability and variability of the stochastic environment that generates proposals for voting. Heterogeneous societies are characterized by the proportions of socially-oriented participants (altruists) \cite{5,6} and other groups of agents. 
However, in this paper we focus only on homogeneous societies.

Previous studies of the ViSE model show that, for all its simplicity, it allows to identify some significant social phenomena and in some cases, to assess their intensity. The analysis of this model helps to answer the following type of questions:
``What happens to a society functioning in accordance with the simplest laws of interaction of human wills under certain conditions?''

We now list the main features of the ViSE model. \emph{Society\/} consists of $n$ \emph{participants\/} (\emph{agents\/}, \emph{voters\/}). Each agent is characterized by the current value of the \emph{capital\/} (a negative capital is interpreted as debt), which changes in discrete time and can also be interpreted as individual utility (an abstract analogue of capital). The vector of participants' capitals before all votes is the initial condition. A \emph{proposal\/} is a vector of proposed \emph{capital increments\/} of all agents. According to the ViSE model, these increments are, in the simplest case, independent and identically distributed random variables, whose mean and standard deviation are denoted by $\mu$ and $\sigma,$ respectively.

The \emph{extinction mode\/} is a version of the ViSE model that presupposes the elimination of the participants whose capitals became negative. In the extinction mode, such participants are called \emph{bankrupts}; in the \emph{no-extinction mode\/} there are no bankrupts, and the agents with negative capitals remain active.

On each step $m\in\mathbb N$:
\begin{itemize}\Up{.6}
\item
The agents who went bankrupt on step $m-1$ are eliminated (``die out'');
\item
If any agents are still active, then the \emph{environment\/} 
generates an independent proposal $\boldsymbol{\zeta}\left(m\right)$; it is a real vector containing one component for each non-bankrupt agent;
\item
The proposal is put to a vote. Each non-bankrupt agent participates in it with one vote;
\item
The proposal is accepted or rejected by means of a certain \emph{voting procedure\/};
\item
If the proposal is accepted, then it is implemented: the capitals of the agents receive the proposed increments.
\end{itemize}

Let $\mathbf{c}(m)$ be the vector of participants' capitals at the end of step $m.$
It has the following dynamics:
\begin{gather}
 \mathbf{c}\left(m\right)=\tilde{\mathbf{c}}\left(m-1\right)+\boldsymbol{\zeta}\left(m\right) I\left(m\right),\quad 0<m\le M,
\end{gather}
where ${\tilde{\mathbf{c}}\left(m-1\right)=\mathbf{c}\left(m-1\right)}$ in the no-extinction mode and $\tilde{\mathbf{c}}\left(m-1\right)$ is obtained from $\mathbf{c}\left(m-1\right)$ by the elimination of all negative components in the extinction mode;
$I(m)$ is $1$ if the proposal $\boldsymbol{\zeta}\left(m\right)$ is accepted and $0$ otherwise; $M$~is the planned length of the ``game.''

The study of the statistical behavior of the capital vector (1) in various environments with various strategies of participants and voting procedures enables one to find the best strategies of the agents and the optimal mechanisms for making collective decisions.

Earlier the competition of individual participants and a group and mechanisms of cooperation in the ViSE model were studied in \cite{1,2,3}; the competition between groups in \cite{3,4}; the influence of altruists on the social dynamics in \cite{5,6}; the ``pit of losses'' paradox (consisting in the fact that a series of democratic decisions may systematically lead the society to general bankruptcy) and dependence of the optimal voting threshold on the environment parameters in \cite{7}; optimization of the group strategy in~\cite{8}.

Relations of the ViSE model to some other models was briefly discussed in~\cite{8}. Here we add the following.
The ViSE model consists of five main elements: 1) repeated voting; 2) the absence of finite ideal states of the agents; 3) individual utilities that enable summation; 4) stochastic environment that generates agenda; 5) social attitudes of the participants (egoism, altruism, group egoism, lobbying, etc.) that determine their voting strategies.

Dynamic models of voting in a multidimensional space of proposals have been intensively studied since the 1960s \cite{9,10,11,12}. However, proposal generation was usually the prerogative of the participants, either of the voters themselves (endogenous agenda): with the choice of a random proposer at every step \cite{13,14,15,16} or not completely randomly \cite{17,18}, or of some other persons with their own interests \cite{9,19}. Proposals generated exogenously and having random effects were considered in \cite{20}\footnote{In \cite{20}, the importance of analyzing models of dynamic voting in a stochastic environment was emphasized and the scarcity of the corresponding papers was noted. The authors mentioned publications since 2008. Note that the ViSE model, which is essentially different from the model of \cite{20}, but belongs to the same class, has been studied since 2003~\cite{21}.} and \cite{22,23}, but the non-transferable nature of the individual utilities (i.e., measuring utilities in subjective scales rather than in a common ``currency'') made it difficult to study cooperative and altruistic strategies within the corresponding models. The studies of cooperation in dynamic voting in the game of pie sharing with transferable utilities started with~\cite{17}. An experimental study of cooperative solutions was conducted in~\cite{24}. Important results on coalition dynamics and decision efficiency in the model with discounting and bribes were obtained in~\cite{15}. The results of \cite{7}, where we looked for the optimal threshold of qualified majority, are comparable with the results of \cite{25}, where a stochastic model was studied for which, under certain conditions, simple majority was optimal. The problem of the optimal voting threshold was also considered in \cite{26}, where proposals were generated exogenously and the profitability of the first accepted proposal was evaluated; in this model, the unanimity procedure is not optimal because of discounting.

An important difference between the ViSE model and legislative bargaining models (most of the papers mentioned above belong to this trend, see also~\cite{27}) is that the question of equilibrium is the main focus of bargaining models, while the structure of proposals' space in the ViSE model, where the coordinates are the individual utilities, does not allow any equilibria. The latter is caused by the lack of Pareto optimal alternatives in the space of opportunities, as the basic ViSE model allows simultaneous gain of all agents in any state. In such circumstances, the goal can be the efficiency of decision making, i.e., maximization of utility: individual, total or average, for certain categories of participants or for the whole society. Using these criteria, it is convenient to study the efficiency of participants' egoism, lobbying, and various altruistic strategies.

To sum up, one of the differences between the ViSE model and the models studied in many other papers is that the focus of research on seeking equilibrium in legislative bargaining is shifted to maximizing some aggregated utility (or capital) in the ViSE model, i.e., to the \emph{efficiency\/} under conditions that can be described in terms typical of games against nature.

The analysis of efficiency brings this approach closer to the one implemented in \cite{28}, where there is no dynamics, since a number of (exogenously generated) proposals are voted ``in parallel,'' but the acceptance of any proposal yields some increments of the individual utilities, and the characteristics of the resulting increment vector make it possible to assess the efficiency of the decision making mechanism.

The paper is organized as follows.
In Section 2, we discuss the problem of analyzing the dynamics of homogeneous societies. In sections~3, 5, and~6, we study the dynamics of the society's capital in the cases where proposals obey normal distributions, symmetrized Pareto distributions, and Student's $t$-distributions, respectively. Section~4 contains technical results regarding the symmetrized Pareto distributions. Section~7 relates to the case of distributions with super-heavy tails; Section~8 examines capital dynamics in a favorable environment. Section~9 summarizes the results of the study, and Section~10 is a brief conclusion.

\section{The problem of homogeneous societies' dynamics}
\label{s_onEnv}

As noted above, the paper is concerned with homogeneous societies, i.e., societies consisting of egoists or altruists with the same support window. An egoist votes for a proposal if and only if it increases his/her capital. We assume that the decisions of the society are adopted by a simple majority of votes. 

Consider the following altruistic strategies. First, let us order the agents by the increase of the current capital. Since any capital increment is a realization of a random variable with a continuous distribution, the current values of capital can be assumed to be different, provided that there was at least one adopted proposal. Suppose that a positive integer ${n_0\le n}$ is fixed. For a current proposal, find the sum of the capital increments of $n_0$ poorest agents. An altruist supports the proposal if and only if this sum is positive, regardless of the capital changes of the remaining ${n-n_0}$ participants. In the case of ${n_0=n},$ the altruist supports exactly the proposals that increase the total capital of the society; for a smaller $n_0,$ the proposals that enrich, in total, the poorest strata consisting of $n_0$ agents. If the initial capital of the agents was the same and no proposals have been adopted yet, the agents are technically ordered by their number and this order is used in the altruistic strategy with threshold~$n_0$.

Let the initial capital of each participant be $C_0$;
$\mu$ and $\sigma$ are the mean and standard deviation of the capital increments in the proposals generated by the environment.
As noted above, the model (1) is considered in two versions: the extinction mode and the no-extinction mode.
The ``game'' ends after $M$ voting steps or earlier if all participants went bankrupt.
The \emph{Average one-step Capital Increment\/} (in what follows, \emph{ACI\/}) of the agents is calculated for the whole ``game'' and takes into account not only nonzero increments from the adopted proposals, but also zero increments from the rejected proposals. It does not take into account the zero increments of bankrupts starting from the moment (step) of their elimination. ACI will be the main efficiency criterion of voting strategies.
In the extinction mode, another important criterion is the proportion of participants who did not go bankrupt throughout the game (the \emph{survival rate\/}).

Let us study how the comparative efficiency of altruistic and egoistic strategies depend on the threshold $n_0,$ along with the external parameters $\mu$, $\sigma$, and the type of capital increment distribution. More specifically, it makes sense to find the most efficient values of the threshold $n_0$ for favorable ($\mu>0$), neutral ($\mu=0$), and unfavorable ($\mu<0$) environments and to compare the result with the efficiency of the egoistic strategy.

Solution of such problems can be reduced to finding the mathematical expectation of the ACI and the survival rate under different combinations of conditions, i.e., to multiple integration whose results can rarely be expressed in standard functions. That is why the results presented below are mainly obtained by means of simulation. Note that during the simulation, random samples were generated anew for each experiment and were not  reused.

\section{Dynamics with normally distributed proposals}
\label{s_onThr}

Suppose now that the environment proposals are realizations of independent normally distributed random variables with parameters $\mu$ and~$\sigma$.

\begin{figure}[t]
\centering{\includegraphics[width=0.95\textwidth]{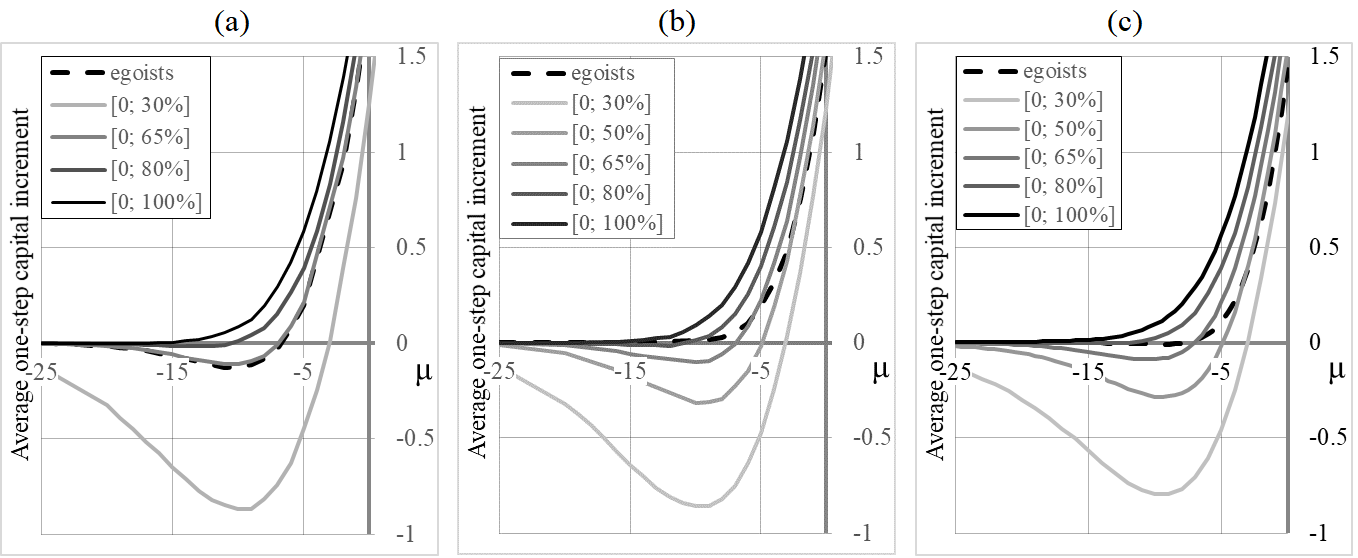}}

    \caption{Average one-step capital increment (ACI) without bankruptcy; number of participants $n=201$; $\sigma=80$. Homogeneous societies consist of egoists or altruists with various support windows. a:~The proposals obey the normal distribution, b:~the symmetrized Pareto distribution with $k=20$, c:~the Student's $t$-distribution with three degrees of freedom.}
\end{figure}

Fig.~1{\rm a} shows ACI (average one-step capital increment) curves (obtained by means of computer simulation aimed at estimation of the corresponding mean) 
in the no-extinction mode.

If the ACI and $\mu$ are measured in $\sigma$ units, 
then the curves in Fig.~1{\rm a} will be invariant to variations of~$\sigma$~\cite{7}. 

For definiteness, let $\sigma=80$; the number of participants in all experiments is $n=201$. If $\mu$ ranges from $-20$ to $-7$, then the ACI in a society of egoists (dotted line) is negative. Thus, the society is getting poorer, although a majority is getting richer on each step!

This phenomenon is known as the ``pit of losses'' paradox; it was studied analytically in~\cite{7}. The cause of the paradox is that the negative increments generated by the unfavorable environment have, on average, larger absolute values than the positive ones. As a result, the total capital increase of the supporting majority (this majority usually slightly exceeds 50\%), on average, is less than the total capital loss of the minority that voted against the proposal. With an increase of the absolute value of $\mu<0$, the percentage of accepted proposals tends to zero and the society of egoists stops to lose capital.

Fig.~1{\rm a} also shows the curves corresponding to the altruistic societies with various support windows. If the support window is narrower than $[0; 80\%]$, then the pit of losses is not negligible. In this case, the total loss of agents who remain outside the support window systematically exceeds the total positive increment of the participants who belong to the support window for a certain~$\mu$.

The ACI curve corresponding to a society of egoists has the same form as the those for altruistic societies, and it is close to the curve for the altruists having support window\footnote{It is even slightly closer to the curve for altruists with support window $[0; 64\%]$; the same holds for $n=801$; when $n=23,$ the closest curve is that for the altruists with the window $[0; 61\%]$.} $[0; 65\%]$. The ACI increases with an increase of the right endpoint of the interval.

\begin{figure}[t]
\centering{\includegraphics[width=0.95\textwidth]{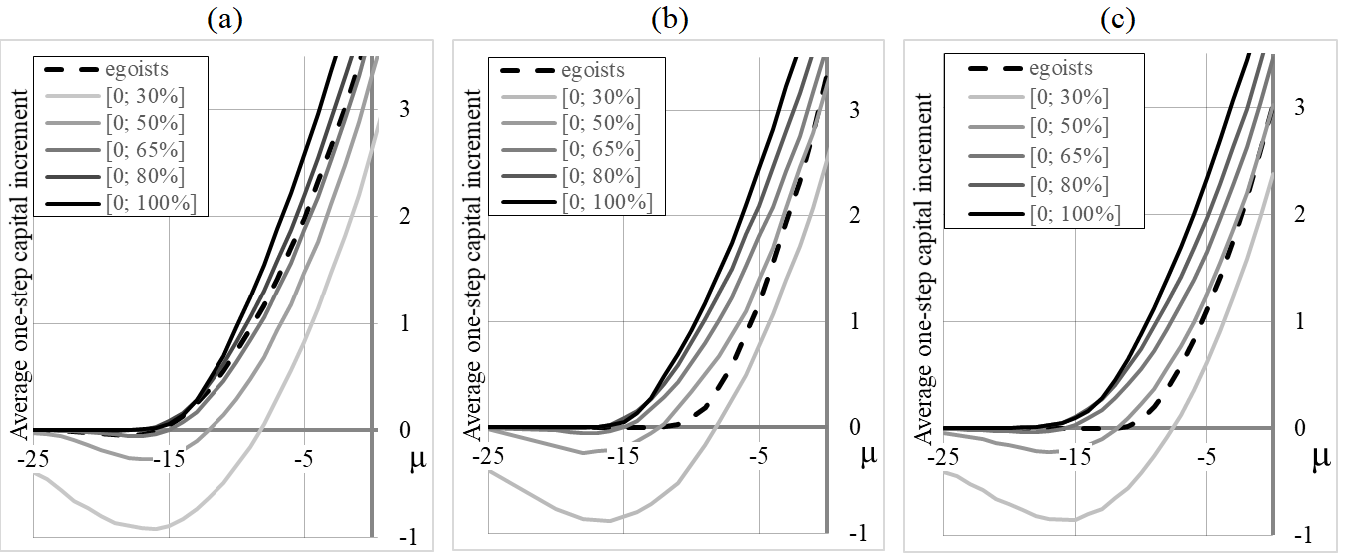}}

    \caption{Average one-step capital increment (ACI) in extinction mode; the number of agents $n=201$; the initial capital of each agent $C_0=40$; $\sigma=80$. Homogeneous societies consist of egoists or altruists. a:~The proposals obey the normal distribution, b:~the symmetrized Pareto distribution with $k=20$, c:~the Student's $t$-distribution with three degrees of freedom.}
\end{figure}

In the extinction mode (the proposals are still generated by the normal distribution), the egoists' ACI curve with the initial capital of each agent $C_0=40$ is also quite close to the curve of altruists with support window $[0; 65\%]$ (Fig.~2{\rm a}), however, the egoists' curve goes higher.

\begin{figure}[t]
\centering{\includegraphics[width=0.95\textwidth]{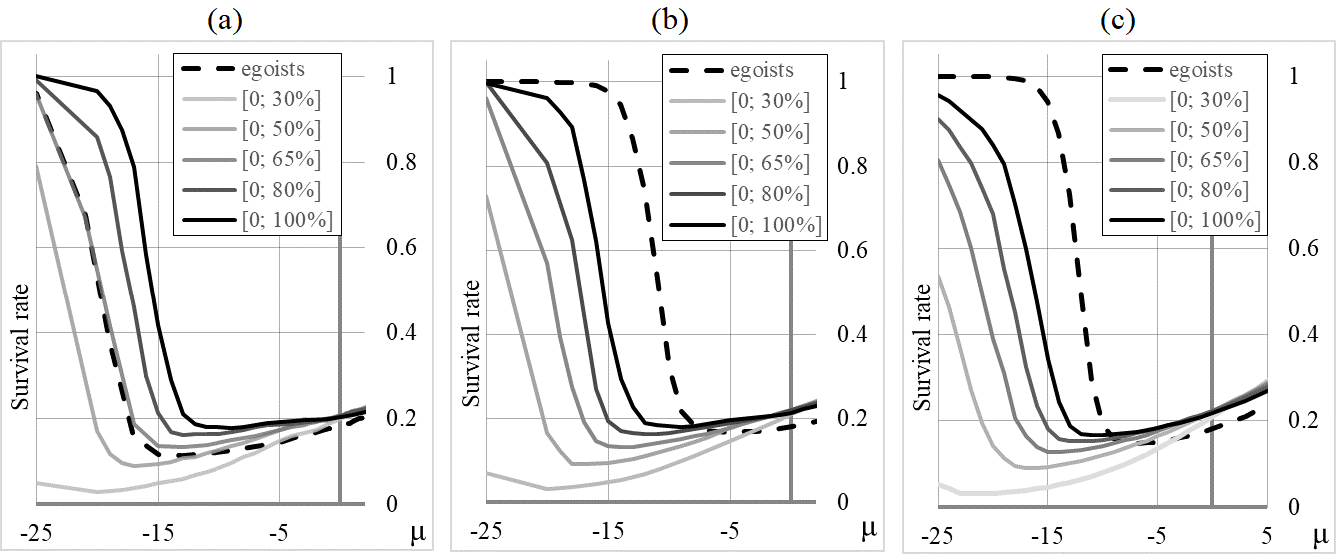}}

   \caption{The proportion of agents at the end of the game (500 steps) in the extinction mode; the number of agents $n=201$; the initial capital of each agent $C_0=40$; $\sigma=80$. Homogeneous societies consist of egoists or altruists. a:~The proposals obey the normal distribution, b:~the symmetrized Pareto distribution with $k=20$, c:~the Student's $t$-distribution with three degrees of freedom.}
\end{figure}

The curves of these two societies are also quite similar on the diagram showing the dependence of the survival rate on the mean $\mu$ of the environment proposals (Fig.~3{\rm a}). At the same time, the results of egoists are slightly worse when $-20<\mu<0$. 

Are these results general, are they preserved for the other distributions of environment proposals? To answer this question, in Section~4 we consider a family of distributions that are substantially different from the normal distributions, namely the symmetrized Pareto distributions.

\section{Symmetrized Pareto distributions}
\label{s_SymPare} 

Pareto distributions are distributions with heavy right tails\footnote{Normally, this concept refers to the distributions such that
$\int^{\infty }_{-\infty }{e^{tx}dF\left(x\right)=\infty}$ for all $t>0$, where $F\left(x\right)$ is the cumulative distribution function.}, which are in use since Pareto's works (see \cite{29}) for econometrics (modeling welfare and income distributions, prices of financial assets, etc.), linguistics, etc.

If a random variable $X$ has a Pareto distribution with parameters $k>0$ and $a>0$, then $X\in[a,\infty)$ and $P\left(X>x\right)={\left({{\frac{x}{a}}}\right)}^{-k}.$ The mean and variance of $X$ are: ${\mu=\frac{ka}{k-1}}$ (provided that $k>1$) and ${\sigma^2=\frac{\mu^2}{k(k-2)}}$ (provided that $k>2$), respectively.

Pareto distributions are characterized by a uniform relative variability in the domains of high and moderate values. This feature is called scalability or self-similarity and determines the scope of applications of these distributions. Another important characteristic of them is their ``heavy'' tails:  the density decreases, as $x \to \infty$, following the power law, i.e., essentially slower than the densities of exponential or normal distributions.

Realizations of Pareto-distributed random variables are always positive, while the ViSE model assumes that the capital increments can be both positive and negative. Furthermore, distributions that are symmetric with respect to their modes are easily comparable with the normal distributions. That is why we will use \emph{symmetrized\/} Pareto distributions.

This class is constructed as follows. Let us extend the density function $f(x)=\frac{k}{x} \left(\frac{a}{x}\right)^k,$ $x\ge a,$ of a Pareto distribution to the whole real axis by combining it with its reflection w.r.t. the line ${x=a}$. Now dividing the resulting function by $2$ gives a symmetric density with the mode ${x=a}.$ Finally, an arbitrary mode $\mu$ can be obtained by a shift yielding
\begin{gather}
\label{e_SP}
f\left(x\right)={\frac{k}{2a}\left(\frac{\left|x-\mu\right|}{a}+1\right)}^{-(k+1)},\quad x\in {\mathbb R}. 
\end{gather}

Indeed, the density of the new distribution at $x$ is equal to the half of the Pareto density at $y$ such that $|x-\mu| = y-a$.
Substituting $y = |x-\mu| + a$ into the Pareto density $f(y)={{\frac{k}{a}}}\ {\left({{\frac{y}{a}}}\right)}^{-(k+1)}$ (${y\ge a}$) and dividing it by $2$ we obtain~(2).

We call the distributions with density (2) the \emph{symmetrized Pareto distributions\/}\footnote{Note that these distributions differ from the ``symmetric Pareto distributions'' considered in \cite{30,21,32} and some other papers and generalize the distributions introduced in \cite{33} (for which ${\mu=a}$). Distributions \eqref{e_SP} were mentioned in \cite{34}, however, in that paper, factor $a^{-1}$ was missing in the corresponding formula on p.~3. Such distributions with ${\mu=0}$ are sometimes called double Pareto distributions, but this term is also used with other meanings.} (\emph{SP-distributions}) with parameters $a$, $k,$ and $\mu$. As well as Pareto distributions, they have a mean (equal to $\mu$ in the case of SP-distributions) whenever $k>1$ and a variance whenever $k>2$. Let us find a relation of this variance to the parameters $a$ and~$k$.

\begin{lemma}
\label{l_1}
The symmetrized Pareto distribution with parameters ${k>2}$ and ${a=\sigma\sqrt{{{\frac{(k-1)(k-2)}{2}}}}}$ has variance $\sigma^2$.
\end{lemma}

The proofs are given in the Appendix.

\begin{figure}[t]
\centering{\includegraphics[width=0.6\textwidth]{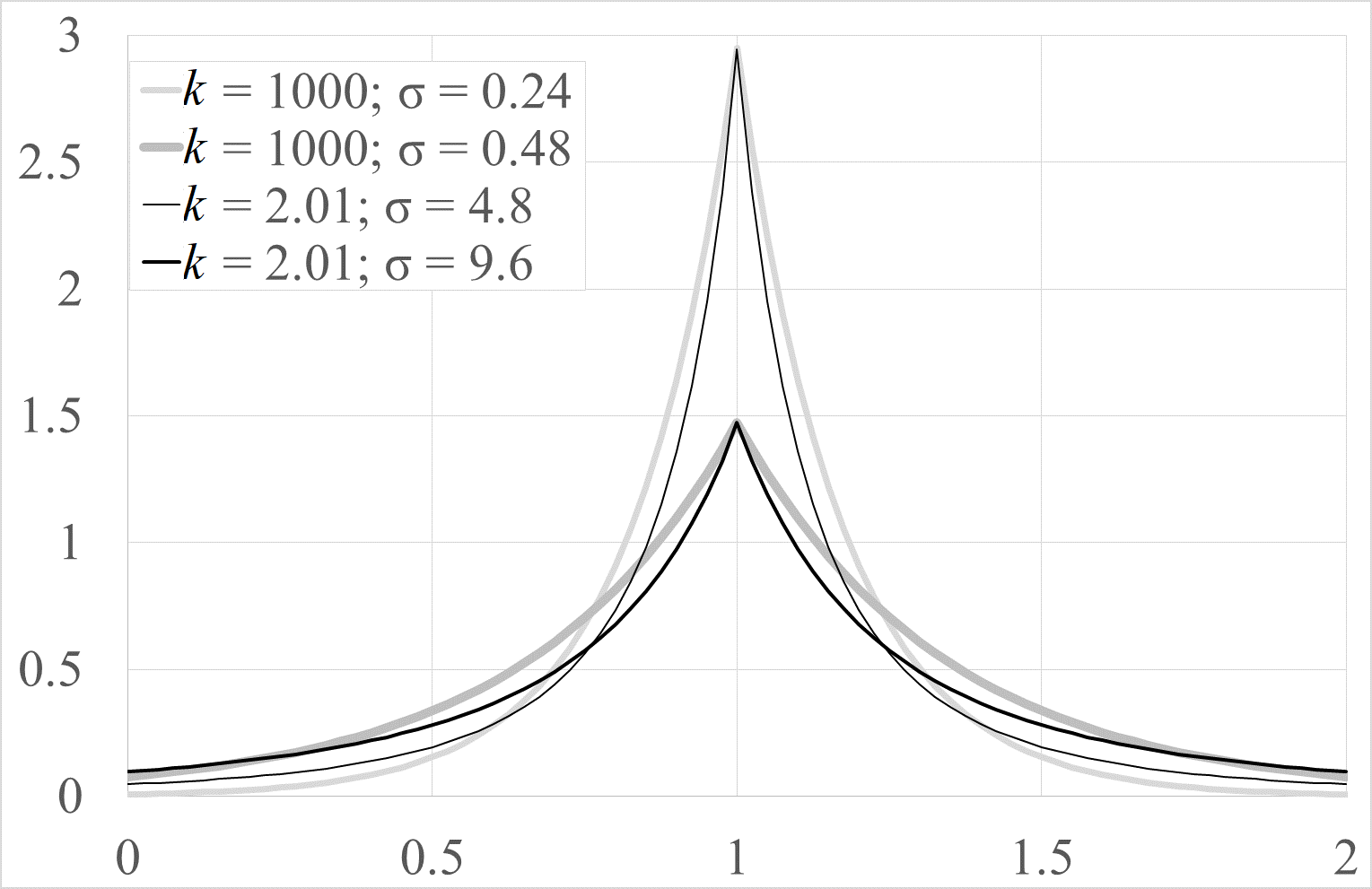}}

     \caption{Examples of density functions of symmetrized Pareto distributions with $\mu =1.$}
\end{figure}

By Lemma~\ref{l_1}, any SP-distribution with $k>2$ can be determined by the parameters $\mu$, $\sigma^2$ and $k$ (the mean, variance, and ``tail lightness,'' respectively). Examples of the density functions of SP-distributions with $\mu = 1$ are given in Fig.~4.
For comparison, two pairs of distributions with the largest density close to $1.5$ and $3$ are presented. The values of $k$ are very high (${k=1000}$) or very low among the distributions having a variance (${k=2.01}$).
For each pair of distributions whose densities are relatively similar in the vicinity of the mode point, the variances differ by a factor of 400. Of course, heavy-tailed distributions have larger variances.

\begin{lemma} 
\label{l_2}
The cumulative distribution function of the SP-distribution with parameters $k>2,$ $\mu,$ and $\sigma$ is
\begin{gather} 
\label{e_SPsi}
F_k\left(x\right)=\left\{ \begin{array}{ll}
\dfrac{1}{2}{\left(1-\dfrac{x-\mu}{\sigma}{{\sqrt{{{\dfrac{2}{(k-1)(k-2)}}}}}}\right)}^{-k}, & x\le \mu, \\[1em]
1-\dfrac{1}{2}{\left(1+\dfrac{x-\mu}{\sigma}{{\sqrt{{{\dfrac{2}{\left(k-1\right)\left(k-2\right)}}}}}}\right)}^{-k}, & x>\mu. \end{array} \right.
\end{gather}
\end{lemma}

To what distribution does the symmetrized Pareto distribution tend as $k \to \infty$? The answer is given by the following lemma.

\begin{lemma} 
\label{l_3}
The SP-distribution with parameters $k> 2,$ $\mu,$ and $\sigma$ tends$,$ as $k\to\infty,$ to the Laplace distribution with the cumulative distribution function
\begin{gather*}
F\left(x\right)=\left\{ \begin{array}{ll}
\dfrac{1}{2}\exp\left({{\dfrac{(x-\mu)\sqrt{2}}{\sigma}}}\right), & x\le \mu, \\[1em]
1-\dfrac{1}{2}\exp\left(-{{\dfrac{(x-\mu)\sqrt{2}}{\sigma}}}\right), & x>\mu. \end{array} \right.\end{gather*}
\end{lemma}

It can be observed that the SP-distributions with $k= 1000$ shown in Fig.~4 are very close to the corresponding Laplace distributions.

The SP-distribution with $k=2$ has no variance. The function $F_k(x)$ defined by \eqref{e_SPsi} converges pointwise (if $k\to 2^+$ and $\sigma$ is fixed) to the function ${\frac{1}{2}\left(\sgn(x-\mu)+1\right)},$ which has a discontinuity at $x=\mu$.

\section{Social dynamics in the case of SP-distribution} 
\label{s_DynSymPar} 

Fig.~1{\rm b} presents the curves that are similar to those shown in Fig.~1{\rm a}, but now the proposals have the symmetrized Pareto distribution. The first thing to note is that for $k=20,$ there is no pit of losses for the society of egoists, and this is a fundamental difference. When $\mu$ increases from $-10$, the ACI of egoists grows slower in the case of an SP-distribution than for a normal distribution.

Fig.~1{\rm b} also shows results for societies consisting of altruists with various support windows. The pit of losses is not negligible when the support window is smaller than $[0; 80\%],$ as well as in the case of normal distribution.

Fig.~2{\rm b} demonstrates the corresponding results for the extinction mode. Compared to the no-extinction mode, the pit of losses (in the societies of altruists) is shifted to the left and  shallower in the case of sufficiently wide support windows. The ACI is noticeably higher for higher $\mu$ than in the no-extinction mode. The ACI of egoists is also higher in the extinction mode. This fact can be explained using the law of large numbers: if the number of participants decreases, then their average capital increment from each adopted proposal exceeds $\mu$ more and more noticeably, and if $\mu<0,$ then the share of accepted proposals grows as well.

Fig.~3{\rm b} presents the curves of the survival rate in the case of extinction. The curve of egoists is dramatically different from that in Fig.~3{\rm a}. With the normal distributions, the survival rate of egoists is slightly lower than that of altruists with support window $[0; 65\%]$, but under the SP-distributions with $\mu <-8$, the egoists' protection from extinction is much better than that of altruists, even when the latter have a support window $[0; 100\%]$.

These results lead to the following conclusion: \emph{egoism better protects society from extinction than altruism in aggressive $($i.e.$,$ substantially unfavorable$)$ environments}. At the same time, with $\mu$ increasing from $-8$, the egoists ``lose'' (by the criterion of survival rate) to altruists with a smaller and smaller support window.

It can be noted that the altruists' results are almost the same in the cases of normal and SP distributions.

The above observations can be explained as follows.
Proposals generated by the SP-distributions with $k=20$ (whose densities decrease rather slowly at a significant distance from $\mu$) have much more outliers (i.e., values distant from the mean) 
than in the samples of the normal distributions with the same~$\sigma$. Each outlier makes a significant contribution to the sample variance, i.e., a smaller part of a fixed variance belongs to the majority of sample values than in the case of the normal distributions. Consequently, most elements of the SP-sample are in a narrow neighborhood of~$\mu$. If $\mu$ is negative and comparable in absolute value with the radius of this neighborhood, then the positivity of most elements of the SP-sample is unlikely, which reduces the chances of a proposal to be accepted. Thus, the \emph{status quo\/} in a society of egoists is sufficiently stable in the case of SP-distributions, and with such values of $\mu$, this society is much more successful in surviving than in environments with normal distributions.

The same factor, i.e., the low proportion of accepted proposals, implies that the ACI of egoists in a moderately unfavorable environment with an SP-distribution is lower than in the case of normal distributions. Comparison of the altruists' curves with support window $[0; 100\%]$ with those of the egoists shows that a large proportion of the proposals that can enrich the society receive insufficient support under the SP-distribution, which leads to rejection of these proposals.

To accept a proposal, altruists need a positive total capital increment of the agents that belong to the support window. It is the sum of a large number of identically distributed random variables, consequently, the curves of altruistic societies under various distributions are identical. For example, for $\mu=-15$ (more precisely, for $\mu/\sigma = 0.1875$, $C_0/\sigma = 0.5$) the survival rate of these societies is quite low and the number of bankrupt agents grows rapidly with the narrowing of the support window. In contrast, egoists successfully survive under the SP-distribution, whereas their average capital remains practically unchanged.

\section{The Case of Student's $t$-Distribution}
\label{s_Student} 

In addition to the heavy tails, the density of the symmetrized Pareto distribution has a sharp peak at the mode and is convex downward at the remaining points (Fig.~4). How do the latter features affect the capital dynamics of the societies consisting of egoists or altruists?

For comparison, we consider distributions with heavy tails and bell-shaped form, namely, Student's $t$-distributions. Those with three degrees of freedom have the heaviest tails among the $t$-distribution with a finite variance. We will call them \emph{S-distributions\/} and consider them as the generator of environment proposals.

A comparison of the tail heaviness of the SP- and t-distributions will be carried out in the next section. The ACI curves of an S-distribution in a no-extinction mode are presented in Fig.~1{\rm c}. The curves of altruists are basically the same as the corresponding curves under the normal distributions and the SP-distributions. With a decreasing $\mu,$ the egoists' ACI decreases faster than under the SP-distributions; the pit of losses is negligibly small.

In the case of extinction, the average capital increment under the S-distributions is shown in Fig.~2{\rm c}. The dependencies on $\mu$ are, in general, the same as for the SP-distribution (Fig.~2{\rm b}). In contrast to the no-extinction mode (Fig.~1), with a decrease of $\mu$ from $-5$ to $-12,$ the ACI under the S-distributions decreases slightly slower than under the SP-distributions; there is no pit of losses as before. When $-3<\mu<0,$ the ACI curve of the egoists under the S-distributions (both with extinction and without it) is closer to the curve of the altruists with support window $[0; 50\%]$ (rather than $[0; 65\%],$ as under the normal distributions).

The curves of survival rate shown in Fig.~3{\rm c} are also similar to the corresponding curves under the SP-distributions (Fig.~3{\rm b}), and the conclusion is the same: the society of egoists is much better protected from extinction in aggressive environments than the society of altruists, however, the former starts losing in a moderately unfavorable environment. In a favorable environment, as we will see in Section~\ref{s_DynamFavor}, 
altruists have a noticeable advantage over egoists. Moreover, low support thresholds are more efficient in terms of survival than high ones.

\section{Social Dynamics with Super-heavy-tailed distributions}
\label{s_Super-heavy} 

Of special interest is social dynamics in the environments whose proposals have super-heavy-tailed distributions. In this study, these are the SP-distributions with $k$ slightly exceeding~$2.$ In this connection, let us recall that an SP-distribution has a variance whenever $k>2$; otherwise the corresponding integral diverges.

\begin{figure}[t]
\centering{\includegraphics[width=0.95\textwidth]{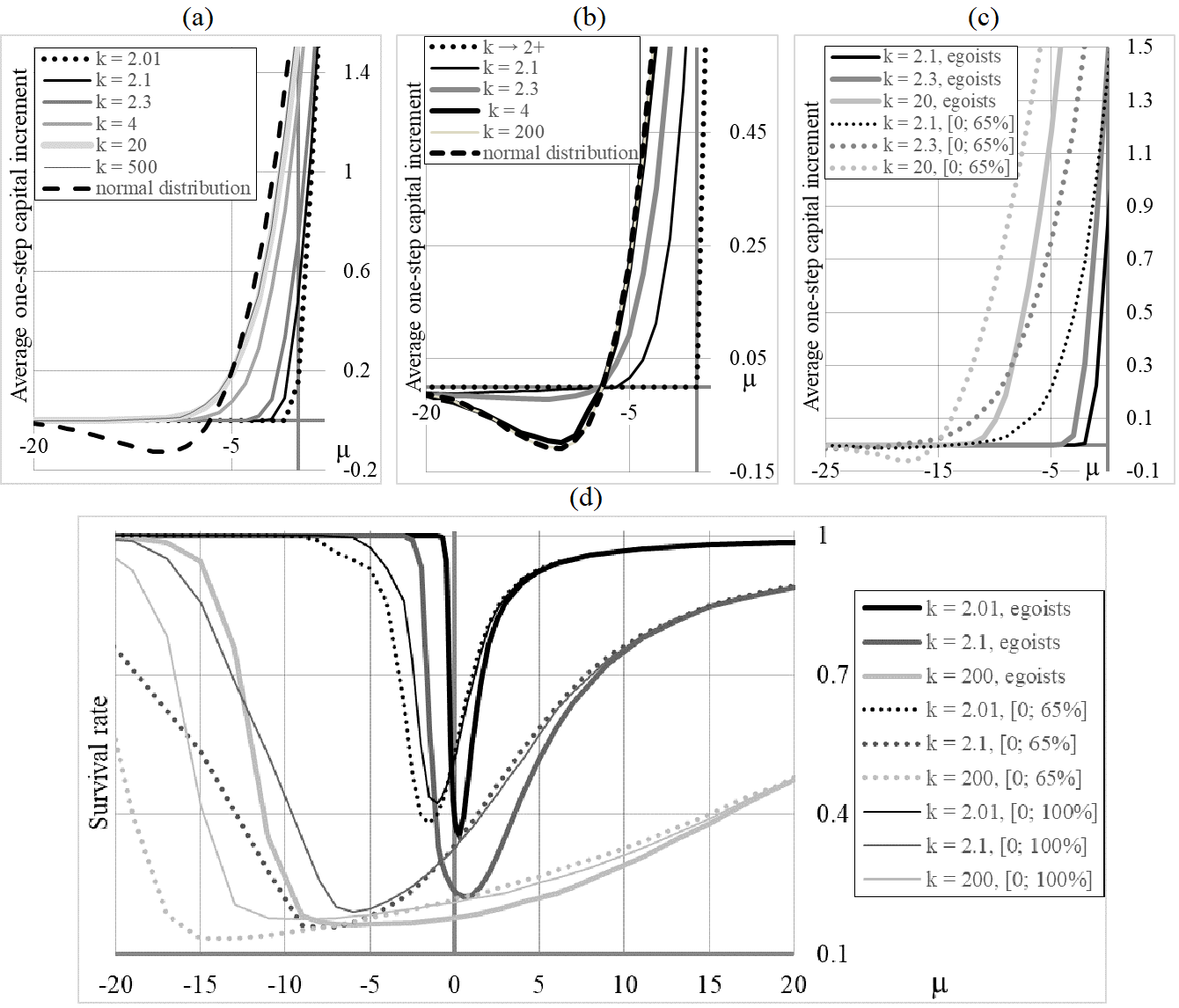}}

    \caption{SP-distributions with $k\in (2; 2.3]$ in comparison with other distributions. The ACI; the no-extinction mode, a:~egoists, b:~altruists with support window $[0, 65\%]$. Egoists and altruists in extinction mode, c:~ACI, and d:~survival rate after 500 steps. The number of participants is $n=201$, the participant's initial capital is $C_0=40$; $\sigma=80$.}
\end{figure}

The results for the distributions with super-heavy tails in comparison with other distributions are shown in Fig.~5.

The curves shown in Fig.~5 correspond to the societies of egoists in the no-extinction mode under the normal distributions and the SP-distributions with various values of~$k$. Note that the curves of egoists under the normal distributions, SP-distributions with $k=20,$ and S-distributions with three degrees of freedom have already been presented in Fig.~1. There is no pit of losses under the SP-distributions with any $k$; the curves for $k=500$ and $k=20$ are close to each other; if $k$ decreases, then the ACI curve shifts to the right and with $k=2.01$ and $\mu<0,$ the proposals are accepted very rarely. Thus, in the case of super-heavy-tailed distributions and negative $\mu,$ egoists protect themselves from extinction very efficiently due to a very rare acceptance of proposals, but they do not increase their capital.

Fig.~5{\rm b} shows the curves of altruists with support window $[0; 65\%]$ without extinction. Earlier, analogous graphs for the normal distributions, SP-distributions with $k=20,$ and S-distributions with three degrees of freedom were presented in Fig.~1; they are very similar; an explanation of this fact was given at the end of Section~5. Basically, this is also true for the general SP-distributions, except for those with low values of~$k$. SP-distributions with super-heavy tails are so ``marginal'' that even the sum of a large number of treir realizations have specific characteristics.

Comparing Figures~5{\rm a} and 5{\rm b}, it can be noticed that altruists with window $[0; 65\%]$ increase their capital faster than egoists if the proposals obey the SP-distributions with various $k$ and $\mu>-7$ ($\mu/\sigma>-7/80$), but with $\mu<-7,$ altruists' survival rate is worse, while egoists retain the \emph{status quo}.

Let us turn to the extinction mode (Figures~5{\rm c} and 5{\rm d}). Here, the ACI of the egoists is nonnegative as well. The ACI curve of altruists with support window $[0; 65\%]$ has a noticeable pit of losses with $k=20$ and this pit disappears with the decrease of~$k$. However, the smaller $k,$ the slower the ACI grows with the increase of~$\mu.$ Furthermore, the ACI of altruists with window $[0; 65\%]$ grows much faster than that of egoists with the same~$k$. Thus, egoism effectively solves the problem of preserving capital in an aggressive environment, however, it is less efficient that altruism (with support window $[0, 65\%]$) in the capital increase of survivors in a moderately unfavorable environment. Note that in the case of distributions with super-heavy tails, the capital preservation problem is irrelevant, since the altruists with window $[0; 65\%]$ have a negligible pit of losses.

In the case of extinction, an essential criterion is survival rate, i.e., the relative number of participants at the end of the game (Fig.~5{\rm d}). It is interesting to note that for super-heavy-tailed SP-distributions ($k=2.1$; $k=2.01$), environments with \emph{small positive\/} values of $\mu$ are the most devastating for egoists. This can be explained by the fact that for \emph{negative\/} values of $\mu$, even with small absolute values, the proportion of accepted proposals is much lower than for positive ones. With $k=200$, the minimum survival rate is reached with a \emph{negative\/} $\mu$, however, a \emph{decrease\/} of $\mu$ from the minimum point gives a larger increase in the survival rate than the same \emph{increase\/} of~$\mu$. At the same time, the higher $k$, the wider an analog of the pit of losses on the survival curve.\x{}

Comparing the survival rate of the egoists and the altruists with support windows $[0; 100\%]$ and $[0; 65\%]$, one can observe that the ``pit'' on the curve of the first category of agents is wider and is located to the left of that of the second one, while for the third category, it is wider because of the left side and deeper than for the second one. With an increase of $k,$ the pit is usually getting deeper. The support window $[0; 65\%]$ better than $[0; 100\%]$ protects the society from extinction only in a part of the right slope of the pit; for the remaining values of $\mu,$ it causes a worse survival rate.

To clarify the nature of the revealed regularities, we now plot a graph of ``tail heaviness'' of the considered distributions. For any random variable that has a variance, consider the function $w(z)$ that assigns to a non-negative argument $z$ the probability that the variable differs from its mean $\mu$ by at least 
$z\sigma,$ where $\sigma$ is the standard deviation of the variable. More formally,
\begin{gather} 
\label{e_tailH}
    w(z) = F(\mu-z\sigma) + 1 - F(\mu+z\sigma),
\end{gather}
where $F(\cdot)$ is the cumulative distribution function of the random variable. 
Note that by Chebyshev's inequality, $w(z)\le z^{-2}.$

Substituting the expression for $F_k(x)$ given by Lemma~2 into (4), we obtain an expression for the tail heaviness function of the SP-distributions. By construction, it does not depend on the parameters $\mu$ or~$\sigma$.
\begin{corollary} 
For the symmetrized Pareto distribution with $k> 2,$ the tail heaviness function \eqref{e_tailH} is
$w\left(z\right)=\left(1+z\sqrt{\frac2{(k-1)(k-2)}}\right)^{-k}.$
\end{corollary}

Substituting the expression for $F(x)$ given by Lemma~3 into (4), we obtain the limit tail heaviness function of the symmetrized Pareto distributions.
\begin{corollary} 
For the SP-distributions$,$ the asymptotics of the tail heaviness function as $k\to\infty$ is $w\left(z\right)=e^{-z\sqrt{2}}$.
\end{corollary}

Note that the SP-distributions with $k\to\infty$ still have rather heavy tails in comparison with normal distributions. Indeed, the density of the Laplace distribution to which the SP-distributions with $\mu=0$ and $\sigma=1$ tend, at a considerable distance from $0$ decrease much slower than that of a normal distribution.

\begin{figure}[t]
\centering{\includegraphics[width=0.7\textwidth]{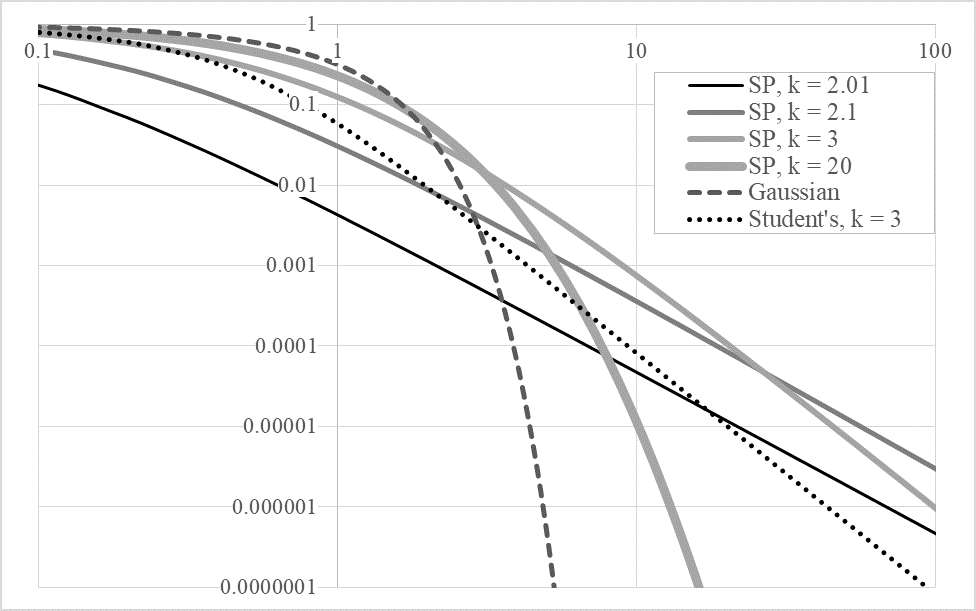}}

    \caption{Tail heaviness functions of the distributions used as proposal generators: symmetrized Pareto distributions (SP), Student's $t$-distribution with three degrees of freedom, and normal (Gaussian) distribution.}
\end{figure}

The tail heaviness functions of the considered distributions in double logarithmic coordinates are shown in Fig.~6.

A close vicinity of 
$\mu$ is the ``lightest'' for the normal distribution, as this distribution is the owner of the heaviest complementing tails.
However, starting with $1.74\sigma$, the tails of the normal distribution successively become lighter than those of the other considered distributions.

On the other hand, the probability of the SP-distributions with $k=2.01$ is mainly concentrated in a small neighborhood of~$\mu$. Say, its realization belongs to the interval ($\mu-0.18\sigma; \mu + 0.18\sigma$) with probability $0.92$. Hence $w(0.18)=0.08$. The SP-distribution with $k=20$ has the ``heaviest'' tails (among the distributions under study) in the interval from $z=1.74$ to $z=3.0$, whereas the SP-distribution with $k=3$ has the heaviest tails from $z=3.1$ to $z=26.61$.

It may seem that this contradicts the view that SP-distributions with $k=2.01$ and $k=2.1$ have super-heavy tails. In fact, this view is adequate for the tails very distant from~$\mu$. For example, the SP-distribution with $k=2.01$ ``overtakes'' the SP-distribution with $k=2.1$ (leading from $z=26.62$) by the tail heaviness only starting with $z=9\cdot10^{10}$,\x{z} while the total weight of its two tails determined by this point is only $5\cdot 10^{-25}$. If the number of voters and the duration of the game are moderate, then sample elements belonging to such tails will never appear in most realizations of the game. Therefore, the main feature of the SP-distribution with $k=2.01$ from the point of view of the ViSE model is that it is largely concentrated in a narrow (in the units of $\sigma$) neighborhood of~$\mu$. Say, its tails are heavier than those of the normal distribution only for $z>3.6$\x{z} (when the tail heaviness is at most $3.6\cdot 10^{-4}$), and are heavier than those of the SP-distribution with $k=20$ starting with $z=8$\x{z} (where the tail heaviness is $7.4\cdot 10^{-5}$).

In a weakened form, these conclusions can be applied to the SP-distribution with $k=2.1$ as well. Therefore, in the context of the ViSE model, it makes sense to consider the SP-distributions with $k=20$ and $k=3$ as the ones with heavy tails. Student's $t$-distribution with three degrees of freedom (overtaking the normal distribution starting with $z=3$\x{z} and the tail heaviness $0.003$) can also be attributed to this class. 
However, it is jointly majorized in the tail heaviness by the SP-distributions with $k=20$ and $k=2.1$.

The above analysis makes it possible to conclude that the SP-distributions with $k$ close to $2$ have super-heavy tails at a great distance from $\mu$, and this distance is so rarely realized in experiments that the tail super-heaviness itself has almost no effect on the ViSE model's social dynamics. What really matters, is their concentration near~$\mu.$ More precisely, these distributions can be called distributions \emph{with super-heavy super-distant  tails}. A note clarifying this conclusion will be given at the end of Section~9.

\section{Dynamics in a Favorable Environment}
\label{s_DynamFavor} 

The previous results, except for those presented in Fig.~5{\rm d}, mainly refer to the case of unfavorable environments, in which the problem of preventing bankruptcy of participants is of unquestionable importance. If the environment is favorable, then the financial dynamics of different societies are distinct either. Moreover, it is interesting to find out whether profitable voting strategies remain beneficial when the favorability of the environment varies. The ACI curves in a favorable environment usually require a separate graphical representation because of a wider range along the ordinate axis.

In the case of normal distributions in the no-extinction mode, the ACI curves continue the curves in Fig.~1{\rm a} monotonically; there are no intersections. For $\mu>0$, the ACI increases with an increase of $\mu$ and the right endpoint of the support window; the curve of egoists remains close to the curve of altruists with support window $[0;65\%]$. Thus, the comparative efficiency of egoistic and altruistic strategies does not depend on the environment favorability. We omit the corresponding diagram because of its triviality.

\begin{figure}[t]
\centering{\includegraphics[width=0.95\textwidth]{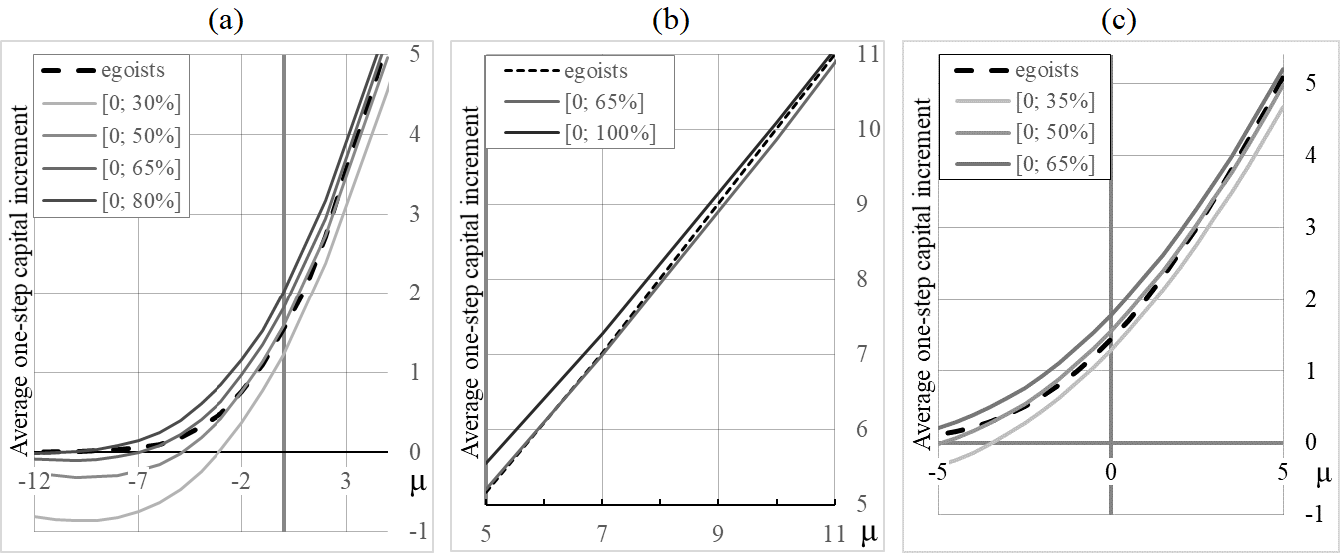}}

     \caption{The ACI without extinction in the case of a favorable environment; the number of participants $n=201$; $\sigma=80$. The proposals obey the SP-distributions with $k=20$ ({\rm a, b}) and Student's $t$-distribution with three degrees of freedom~({\rm c}).}
\end{figure}

For the SP-distribution with $k=20$ (Figures~7{\rm a},{\rm b}), the comparative efficiency of egoistic and altruistic strategies depends on~$\mu$. For $\mu=-12$ the ACI of egoists is higher than that of altruists with support window $[0; 80\%]$, however, for $\mu\in [-2; 2],$ it is lower than even ACI of altruists with support window $[0; 50\%]$. With a further increase in $\mu$, the relative efficiency of egoism increases again and it almost catches up with the efficiency of the altruists with support window $[0; 100\%]$ for $\mu=11$. Thus, we can conclude that under heavy-tailed distributions, egoism is efficient in unfavorable environments, when it is important to avoid individual dangers, and it is also efficient in favorable environments, which regularly 
proposes valuable individual ``prizes.'' But this strategy is too ``fussy'' in environments close to neutral; in this case, even altruism with a rather narrow support window is more efficient. The conclusion is similar for the Student's $t$-distribution (Fig.~7{\rm c}).

\begin{figure}[t]
\centering{\includegraphics[width=0.95\textwidth]{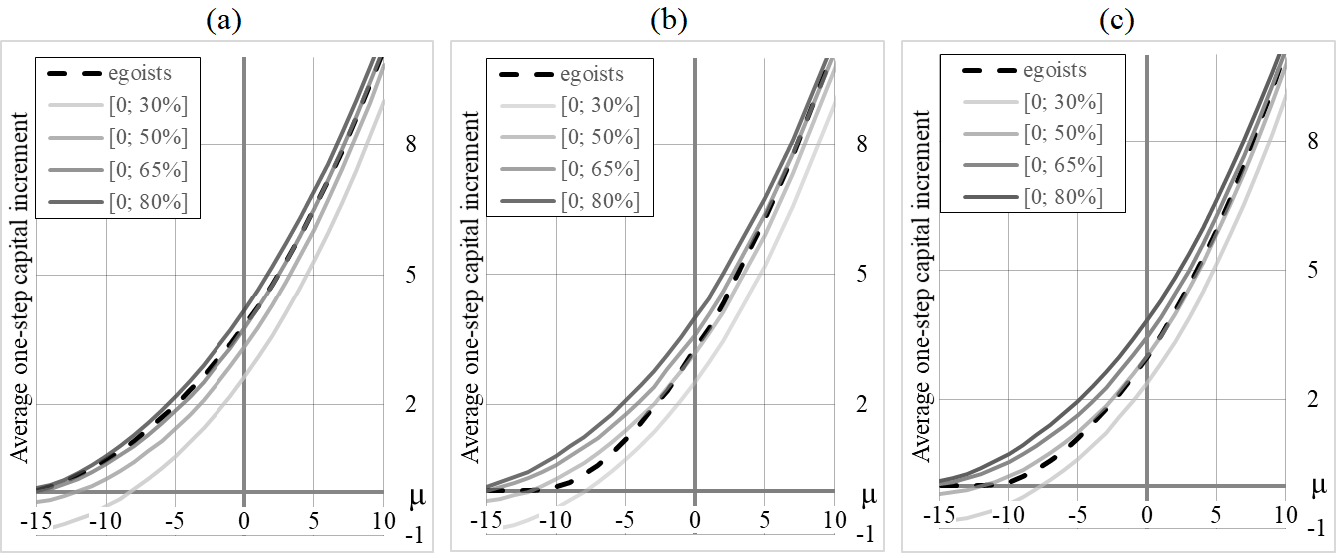}}

     \caption{The ACI in extinction mode in the case of a favorable environment; the number of participants $n=201$; the initial capital $C_0=40$; $\sigma=80$. Proposals obey the normal distribution (a), the SP-distribution with $k=20$ (b), and the S-distribution with three degrees of freedom~(c).}
\end{figure}

Now consider the extinction mode (Fig.~8). With normal distributions, the ACI curve of egoists is close to the curve of altruists with support window $[0; 65\%]$ on all segments of the real axis, however, at $\mu=-5$ it passes a little higher, while at $\mu = 10$ it is slightly lower. For the SP-distribution, the relative efficiency of egoism grows when the environment changes from an unfavorable to a favorable one. With the S-distribution, a similar pattern is observed, but it is less pronounced.

\begin{figure}[t]
\centering{\includegraphics[width=0.95\textwidth]{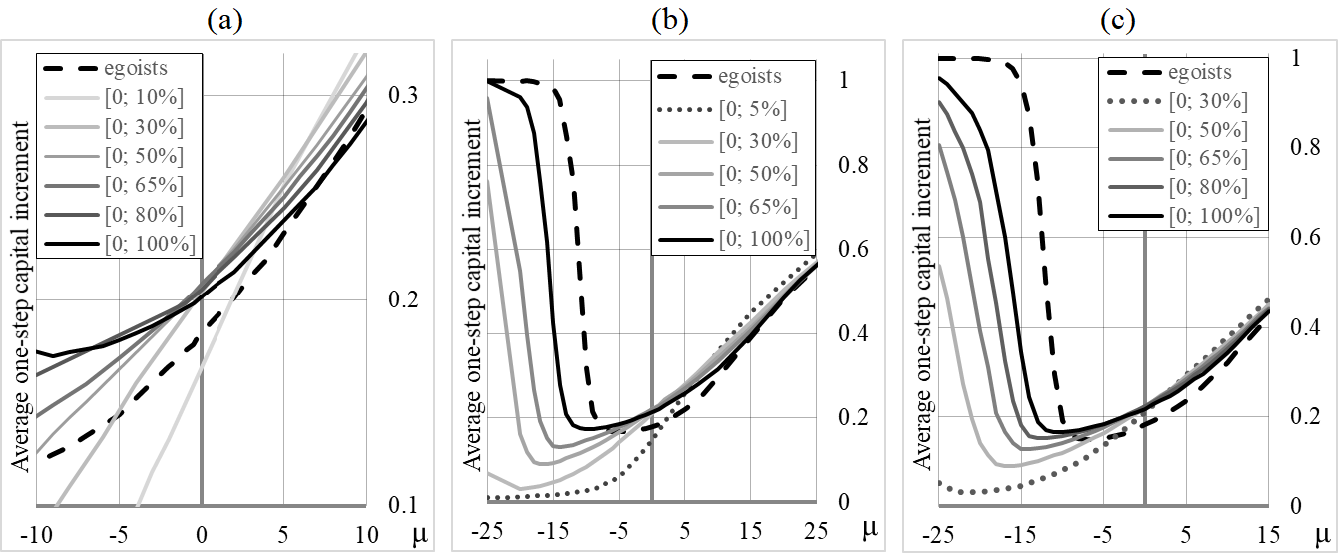}}

     \caption{The survival rate at the end of the game (500 steps): dependence on the environment favorability; the number of participants $n=201$; the initial capital $C_0=40$; $\sigma=80$. Proposals obey the normal distribution (a), the SP-distribution with $k=20$ (b), and the S-distribution with three degrees of freedom~(c).}
\end{figure}

Let us analyze the survival diagram. In the case of normally distributed proposals (Fig.~9{\rm a}), two effects should be noted. First, the compared ACI curves of altruists intersect not far from $\mu=0$, i.e., societies with smaller support window survive more successfully. This confirms the previously formulated principle: \emph{to survive in aggressive environments$,$ it makes sense to support everyone$,$ while in more favorable external conditions$,$ it is better to support only the poorest agents}.

The explanation is simple: in a favorable environment, the ``group of risk'' (i.e., the group of participants with high chances to go bankrupt in one step) is smaller and so it is not efficient to include many agents in the support window, as this would reduce the support of those who are in dire need of it.

Second, \emph{altruism\/} (including altruism with support window $[0; 30\%]$) \emph{is more efficient than egoism in an environment close to a neutral one\/}. Thus, this effect observed in the analysis of Fig.~7{\rm a} is divided into two components in the case of extinction: a less noticeable (ACI, Fig.~8{\rm a}) and a more pronounced (survival rate, Fig.~9{\rm a}).

These two effects are also observed in the cases of SP- and S-distributions (Figures~9{\rm b},{\rm c}). Thus, the effects seem to be robust: they weakly depend on the distribution of proposals.

\begin{figure}[t]
\centering{\includegraphics[width=0.95\textwidth]{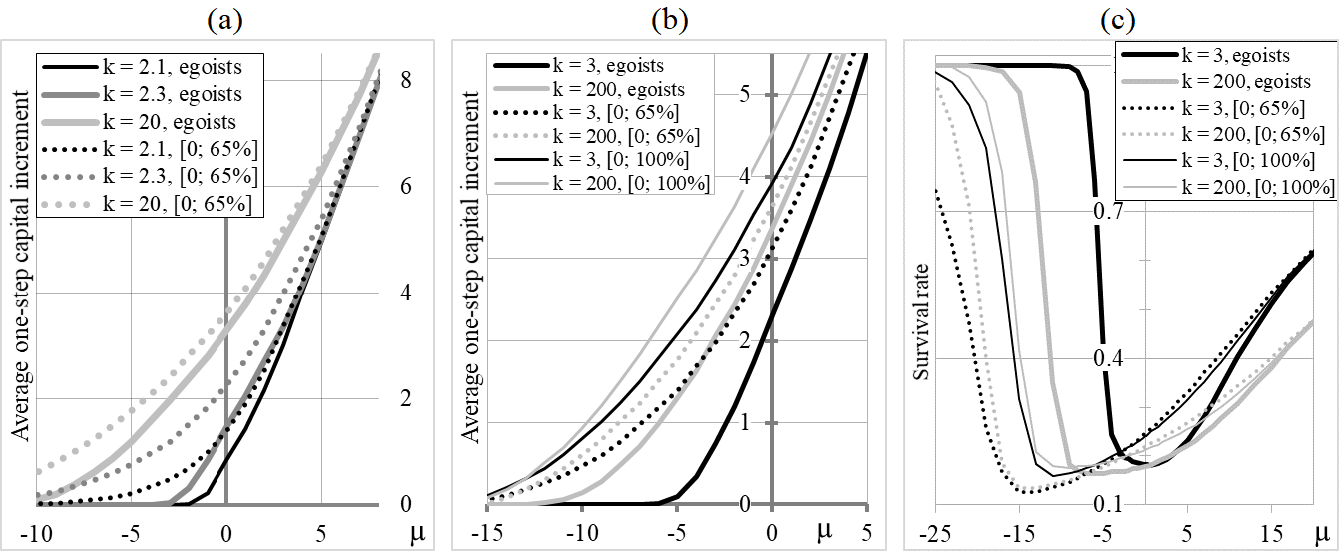}}

     \caption{Indicators of social dynamics with SP-distribution in the extinction mode; $n=201$; the initial capital $C_0=40$; $\sigma=80$. {\rm a$,$ b}: ACI; {\rm c}: the survival rate after 500 steps.}
\end{figure}

Now extend the curves of Fig.~5 to the positive semiaxis (Fig.~10{\rm a}). The curves in Fig.~5  presented the ACI of altruists with support window $[0; 65\%]$ and egoists under the SP-distributions with $k$ close to $2$ and $k=20$ in the extinction mode. For $\mu>0$, similarly as in a moderately unfavorable environment, the average capital grows faster with a higher $k$ and with the altruistic strategy based on the $[0; 65\%]$ window. The support of the whole society (window $[0; 100\%]$) is added to Fig.~10{\rm b} for a high $k=200$ and a rather low $k=3$. This strategy gives the highest ACI. Other results are similar to those in Fig.~10{\rm a}.

Figure~10{\rm c} provides additional material that confirms the conclusions from the analysis of Fig.~5{\rm c}. The following results obtained by comparison Figures~10{\rm b} and 10{\rm c} are  interesting. For all $\mu$ for which the share of accepted proposals is significantly different from $0$ and $1$, the egoists in SP-environments with $k=200$ (relatively light tails for SP-distributions) get a higher ACI than in the environments with $k=3$. At the same time, they have a lower survival rate except for the interval $0.5<\mu<3.8$. Thus, the criteria of ACI and survival rate contradict to each other. In the case of altruism with support window $[0; 65\%]$, these criteria contradict except for the interval $-13<\mu<-7$, where $k=200$ is preferable for both criteria. To the left of this segment, the survival rate is better when $k=200$; to the right, it is better when $k=3$. The situation is similar if the altruistic strategy supports everyone. Generally, SP-distributions with ``lighter'' tails increase the average capital of the egoists faster, while for the average capital of altruists, this is true only for $\mu>-13$. Omitting some details, we can add that with the exception of a certain ``middle'' zone, the type of the distribution tails that ensures a faster growth of the ACI (in different cases, this type varies), simultaneously leads to a faster bankruptcy.

\section{Discussion of the Main Results}
\label{s_overview}

Let us summarize the main results obtained in the research.

1.	On the ``pit of losses'' paradox consisting in the impoverishment of the society in a moderately unfavorable environment due to the implementation of majority decisions. This paradox was originally established for the case of normally distributed proposals.

1.1. There is no pit of losses if the society consists of egoists and the distribution has heavy tails, like the SP-distributions or the S-distributions with three degrees of freedom, see Figures~1{\rm b},{\rm c} and 5{\rm a}.

1.2. With the altruistic strategies under study, a significant pit of losses occurs when the upper bound of the support window is not very high (below 80\%). The shape of the ACI (Average one-step Capital Increment) curves is basically robust to the distribution (Fig.~1), however, the pit of losses becomes shallower in the case of distributions with ``super-heavy tails'' (SP-distribution with $k$ close to~$2$) (Fig.~5{\rm b}).

1.3. In the extinction mode, the above conclusions basically remain valid, the only differences being that there is practically no pit of losses of the egoists' ACI curve in the case of normal distributions, while the pit of losses of the altruists' ACI curve is shifted to the left and becomes shallower when the support window is sufficiently wide (Fig.~2).

2. Comparison of altruistic strategies.

2.1. For the altruists, it is inefficient to support only the poorest in an unfavorable environment: the indices of success are higher if everyone is supported. This conclusion is obvious in the no-extinction mode, since the ACI criterion is maximized by the adoption of all proposals with a positive total capital increment, i.e., by the support window $[0; 100\%]$. However, this conclusion is nontrivial in the extinction mode. It was obtained that narrowing of the support window does not improve well-being (neither the ACI, nor the survival rate) under any distribution considered in this study (Figures 2, 3, 5{\rm d}, and 9). \x{?}

2.2. 
In a favorable environment, smaller support windows better protect from extinction than the larger ones starting with certain values of~$\mu$ (Figures~9 and~10{\rm c}). Nonetheless, large support windows remain more efficient by the ACI criterion in favorable environments too (Figures~8 and~10{\rm b}).

3. Comparison of egoistic and altruistic strategies.

3.1. In the case of normal distribution, the ACI curves of egoists and altruists with support window $[0; 65\%]$ are quite close to each other both in the extinction and no-extinction modes (Figures~1{\rm a}, 2{\rm a}, and 8{\rm a}). In the extinction mode, the first curve slightly lags behind the second one with an increase of~$\mu.$ The curves of survival rate behave similarly, but the egoists are most lagging behind the altruists with support window $[0; 65\%]$ when $\mu$ takes small negative values (Figures 3{\rm a} and 9{\rm a}).

3.2. In the case of normally distributed proposals, altruists with support window wider than $[0; 65\%]$ are more efficient, by both criteria, than egoists, except for the survival rate in a highly favorable environment: in the latter case, their efficiency is minimal, while the leaders (leaving the egoists behind) are altruists with narrow support windows (see the above item 2.2).\x{simplify?}

4. On heavy-tailed distributions: for them, the conclusions are significantly different.

4.1. In the case of heavy tails, egoism is more efficient than altruism in surviving in highly unfavorable environments (Figures~3{\rm b},{\rm c}). By survival rate, egoists are far ahead of even the altruists supporting the entire society. However, they become less efficient in less aggressive and favorable environments.

4.2. In the case of no-extinction mode and heavy tails (Figures 1{\rm b},{\rm c} and 7), the efficiency of egoism (by ACI) in comparison with altruism with various support windows decreases with an increase of $\mu<0.$ Then it reaches the smallest (compared to altruism) value near $\mu=0$ and hereafter grows with a growth of~$\mu$. Thus, \emph{egoism is most justified in highly unfavorable or highly favorable environments$,$ where it occasionally allows agents to avoid individual danger or not to miss personal gain$,$ respectively}.

4.3. With the ACI criterion in the extinction mode, the  results are basically similar (Figures~2{\rm b},{\rm c} and 8{\rm b},{\rm c}), except that egoism is least efficient not near $\mu=0$, but in aggressive environments for which the ACI of egoists takes near-zero positive values. In such environments, unlike more aggressive ones, egoists do not have a significant advantage in survival over the altruists that support everyone (Figures~3{\rm b},{\rm c}).

5. On ``super-heavy-tailed'' distributions (SP-distributions with $k$ slightly exceeding $2$): they produce a specific dynamics.

5.1. In the case of super-heavy tails and a noticeable negative $\mu$, the egoists' ACI is close to zero, as the proposals are adopted very rarely (Figures~5{\rm a},{\rm c}). The ACI of altruists with support window $[0; 65\%]$ is smaller (in absolute values) than the common values for the other distributions under study (Fig.~5{\rm b}), so the pit of losses becomes shallower. In more favorable environments, the ACI of altruists considerably exceeds (except for the most favorable conditions) the ACI of egoists (Figures~5{\rm c} and 10{\rm a}).

5.2. With an increase of tail heaviness, the survival rate of both egoists and altruists becomes higher; the ``pits'' on the corresponding curves become narrower (Fig.~5{\rm d}). It is interesting that among super-heavy-tailed environments, the most ``confiscatory'' for egoists are weakly-\emph{favorable\/} environments.

5.3. The survival curves of altruists supporting the whole society in comparison with those of the egoists are shifted to the left and ``stretched out.'' That is, 
egoism better than altruism protects from extinction in an unfavorable environment, while altruism is more efficient in the environments that provide egoists with a positive ACI; it is the same conclusion as for the other distributions.

5.4. With a decrease of the support window from $[0; 100\%]$ to $[0; 65\%]$, the left ``pit slope'' of the survival curve shifts to the left and the pit becomes deeper, i.e., the latter altruistic strategy is less efficient in unfavorable environments. However, in favorable conditions, as well as for the other distributions, it is somewhat better (Fig.~5{\rm d}).

Finally we make two remarks that help better understand the results.

The first one is about the pit of losses paradox. There is the \emph{algorithm of small handouts to the majority\/} used in the proof of the voting paradox by A.V.~Malishevsky \cite[Chapter~2, Section~1.3]{9}. It bankrupts a society of egoists in a series of voting steps by proposing small handouts to a majority in combination with the ``robbery'' of a minority; each voter is in the minority at one step. The pit of losses paradox (Fig.~1{\rm a}, see also \cite{7}) means that the same algorithm is actually implemented by an unfavorable stochastic environment that generates normally distributed proposals. However, if proposals obey a distribution with heavy tails (such as the SP-distribution or the $t$-distribution with three degrees of freedom, Figures~1{\rm b},{\rm c}), then 
the convergence rate of the average of a large number of realizations to a negative $\mu$ is insufficient to ensure that the loss of a ``wandering'' minority systematically exceeds the gain of the majority that supports adopted proposals. Therefore, the agents survive (or bankrupt negligibly slowly if the pit of losses is very shallow). Paradoxically at first glance, they are protected by outliers in the proposals, sometimes called (after J.St. Mill) ``black swans''~\cite{35}. This effect can be interpreted as a specific manifestation of ``antifragility'' discussed in~\cite{36}.

The second remark concerns distributions with super-heavy tails. The present paper analyzes the ViSE originated social dynamics with distributions of different types, but the same mean and variance. This standardization of parameters is rather common and seems to be natural. However, let us return to Fig.~4 showing two pairs of SP-distributions which are, in each pair, relatively close in a neighborhood of~$\mu$. As was noted, the standard deviations $\sigma$ of these ``roughly similar'' distributions differ by a factor of~20. If we transform the distributions with super-heavy tails ($k=2.01$) in such a way that they acquire the same variance as the distributions with $k=1000$, then the former distributions contract 20 times and lose all resemblance to the latter distributions in the neighborhood of~$\mu$. 
The question is: which distribution with $k=2.01$ should be considered to be similar to the corresponding distribution with $k=1000$? The initial one behaves rather similarly in the central zone, but its variance is 400 times larger. The transformed one has the same mean and variance, but its behavior in the neighborhood of $\mu$ is completely different. There are strong arguments for the first answer. This fact indicates the desirability of an alternative approach to the standardization of distributions. Say, one can consider two symmetric distributions to be ``similar'' when they have the same $\mu$ and the same interval (centered at $\mu$) containing a certain essential proportion of probability. In this case, the tails 
cut off by this interval will have the same heaviness, with a different heaviness of more distant tails. With this approach, the compared distributions will have more similarity in the central region, the super-heavy tails will not be ``super-outlying,'' and the conclusions of the present research will be complemented with new ones.

\section{Conclusion}
\label{s_Concl}

In this paper, we studied the social dynamics determined by the ViSE model and compared the efficiency of egoistic and altruistic voting strategies. We found that the ``pit of losses'' paradox described earlier does not appear under several heavy-tailed distributions. Furthermore, for such distributions, the egoistic strategy better protects society from extinction in unfavorable environments than the altruistic ones. On the other hand, in more favorable environments, the efficiency of altruism is higher than that of egoism. The comparison of altruistic strategies with each other shows that for the maximum survival rate in aggressive environments, altruists must support everyone, and this conclusion is robust to distributions, while in more favorable environments, it is better to support only the weakest agents. It turns out that under heavy-tailed distributions, the comparative efficiency of egoism is, in some cases, higher in markedly unfavorable or markedly favorable environments (where it is important to avoid personal dangers or not to miss individual ``prizes,'' respectively) and it is significantly lower in intermediate environments. In the next step, we are going to consider agents with combined voting strategies and strategically heterogeneous societies.

\section*{Acknowledgments}
The work of P.Yu.~Chebotarev and V.A.~Malyshev was supported by the Russian Science Foundation, project 16-11-00063 granted to IRE~RAS.

\appendix{}

%
%

\PLE{\ref{l_1}}
Since the variance of the symmetrized Pareto distribution for fixed $a$ and $k$ is invariant to the variation of $\mu$, let us find it for $\mu = a$:
\begin{gather*}
\sigma^2=\int\limits^a_{-\infty }{{\left(x-a\right)}^2{{\frac{k}{2a}}}}{\left(2-{{\frac{x}{a}}}\right)}^{-k-1}dx+\int\limits^{\infty
}_a{{\left(x-a\right)}^2{{\frac{k}{2a}}}}{\left({{\frac{x}{a}}}\right)}^{-k-1}dx{}
\\
{} ={{\frac{ka^k}{2}}}\left(\int\limits^{\infty }_a{{\left(a-y\right)}^2}y^{-k-1}dy+\int\limits^{\infty }_a{{\left(x-a\right)}^2}x^{-k-1}dx\right){}
\\
{} =ka^k\int\limits^{\infty
}_a\left(x^{-k+1}-2ax^{-k}+a^2x^{-k}\right)dx{}\\
{}=ka^2\left({{\frac{1}{k-2}}}-{{\frac{2}{k-1}}}+{{\frac{1}{k}}}\right)={{\frac{2a^2}{\left(k-1\right)\left(k-2\right)}.
}}\end{gather*}
Lemma~\ref{l_1} is proved.

\PLE{\ref{l_2}}
Using formula (2) for the density of the SP-distribution and Lemma~\ref{l_1}, for $k>2$, $\mu = 0$ and $x\le 0$, we obtain
\begin{gather*}
F_k\left(x\right)=\frac{k}{2a}\int\limits^x_{-\infty }{\left(1-\frac{y}{a}\right)}^{-\left(k+1\right)}dy{}
{}=\frac{1}{2}{\left(1-\frac{x}{a}\right)}^{-k}
=\frac{1}{2}{\left(1-\frac{x}{\sigma}{{\sqrt{{{\frac{2}{\left(k-1\right)\left(k-2\right)}}}}}}\right)}^{-k}.
\end{gather*}
The expression for $x>0$ is obtained using the symmetry of the distribution density with respect to $0$; expression for $\mu = 0$ are obtained by replacing $x$ with $x-\mu$.

\PLE{\ref{l_3}}
Using Lemma~\ref{l_2} and the ``second remarkable limit,'' for $x\le\mu$ we find
\begin{gather*} \lim\limits_{k\to \infty }
F_k\left(x\right){}
{}=\frac{1}{2}\lim\limits_{k\to \infty }
\left(1+\left(\frac{\sigma}{\mu-x}\sqrt{\frac{(k-1)(k-2)}{2}}\right)^{-1}\right)^{\frac{\sigma}{\mu-x}\sqrt{\frac{(k-1)(k-2)}{2}} \cdot
\frac{-k(\mu-x)\sqrt{2}}{\sigma\sqrt{(k-1)(k-2)}}}{}
\\
{}=\frac{1}{2}\exp\left(\frac{(x-\mu)\sqrt{2}}{\sigma}\right).
\end{gather*}

The expression for $x>\mu$ is obtained using the symmetry of the distribution. The found limit function equals to the Laplace cumulative distribution function with parameters $\mu$ and $\sigma$.



\end{document}